\documentclass[paper=a4, USenglish, numbers=noenddot]{scrartcl} 
\usepackage{cmap}		
\usepackage[utf8]{inputenc}	
\usepackage[T1]{fontenc}	
\usepackage[USenglish]{babel}	
\usepackage[]{microtype}

\usepackage{amsmath}
\usepackage{amssymb}
\usepackage{graphicx}
\usepackage{mathtools}	

\usepackage{aliascnt}		
\usepackage{amsthm}

\usepackage{tikz}			
\usetikzlibrary{matrix,arrows,patterns,intersections,calc,decorations.pathmorphing}

\usepackage[pdftitle={Notes on the Veech group of the Chamanara surface},pdfauthor={Frank Herrlich, Anja Randecker},pdfborder={0 0 0}]{hyperref}   
\usepackage[figure]{hypcap}		

\newtheoremstyle{break}
 {} 
 {} 
 {\itshape} 
 {} 
 {\bfseries} 
 {} 
 {\newline} 
 {\thmname{#1}\thmnumber{ #2}\thmnote{ (#3)}} 
 
\newtheoremstyle{breakdef}
 {} 
 {} 
 {} 
 {} 
 {\bfseries} 
 {} 
 {\newline} 
 {\thmname{#1}\thmnumber{ #2}\thmnote{ (#3)}} 
 
\newtheoremstyle{remark}
 {} 
 {} 
 {} 
 {} 
 {\itshape} 
 {.} 
 {0.5em} 
 {\thmname{#1}\thmnumber{ #2}\thmnote{ {\normalfont (}#3{\normalfont )}}} 

\theoremstyle{breakdef}

\theoremstyle{remark}

\newaliascnt{lem}{definition}  
\newtheorem{lem}[lem]{Lemma}
\aliascntresetthe{lem}

\theoremstyle{break}

\newaliascnt{prop}{definition}  
\newtheorem{prop}[prop]{Proposition}
\aliascntresetthe{prop}

\DeclareMathOperator{\GL}{GL}
\DeclareMathOperator{\SL}{SL}
\DeclareMathOperator{\PSL}{PSL}

\renewcommand{\Re}{\mathrm{Re}}

\author{Frank Herrlich, Anja Randecker}
\title{Notes on the Veech group \\ of the Chamanara surface}
\date{\today}

\begin{document}

\maketitle

These notes follow and extend the proof of the description of the Veech group in \cite[Theorem 4]{chamanara_04}. The main result is that the Veech group of the Chamanara surface is a non-elementary Fuchsian group of the second kind which is generated by two parabolic elements. Most of the calculations were carried out on the train from Marseille to Karls\-ruhe after the wonderful conference “Dynamics and Geometry in the Teichmüller Space” in July 2015.

\section{Review of definitions}

Let us quickly go through the definitions that we use in this setting.

A \emph{translation surface} $(X, \mathcal{A})$ is a connected two-dimensional manifold $X$ together with a translation structure $\mathcal{A}$ on $X$, i.e.\ a maximal atlas on $X$ so that the transition functions are locally translations.
Note that this is a more general definition compared to what was the standard definition for a long time.

Every translation surface $(X, \mathcal{A})$ can be equipped with a metric that comes from pulling back the Euclidean metric via the charts. The metric completion $\overline{X}$ may contain additional points which are called \emph{singularities}.

There are three kinds of singularities:
A singularity is called \emph{cone angle singularity} if it has a punctured neighborhood in $X$ which is isometric to a finite translation covering of a once-punctured Euclidean disk. It is called \emph{infinite angle singularity} if it has a punctured neighborhood in $X$ which is isometric to an infinite translation covering of a once-punctured Euclidean disk. All other singularities are called \emph{wild}.

We say that a translation surface is \emph{finite} if its metric completion is a compact surface and all singularities are cone angle singularities. This is equivalent to all three standard definitions of translation surfaces, except that the singularities are not included in the translation surface in our setting.

\bigskip

A continuous map $f \colon (X,\mathcal{A}) \to (Y, \mathcal{B})$ of translation surfaces is called \emph{affine} if it is locally affine, i.e.\ for every $x\in X$ there exist charts $(U,\varphi)\in \mathcal{A}$ and $(V,\psi)\in \mathcal{B}$ with $x\in U$ and $f(U)\subseteq V$ so that for every $z\in \varphi(U) \subseteq \mathbb{R}^2$ it is true that
\begin{equation*}
 \left( \psi\circ f \circ \varphi^{-1} \right)(z)= A\cdot z+t \text{ for a fixed } A\in \GL(2,\mathbb{R}) \text{ and a fixed } t\in\mathbb{R}^2 .
\end{equation*}

For an affine map $f$, the matrix $A$ is globally the same for all choices of charts and this matrix is called the \emph{derivative} of~$f$.

For a translation surface $(X,\mathcal{A})$, the \emph{Veech group} $\GL^+(X,\mathcal{A})$ of $(X,\mathcal{A})$ is defined as the group of all derivatives of orientation-preserving homeomorphisms of $X$ that are affine with respect to $\mathcal{A}$.

For a translation surface with finite area, the elements of the affine group are area-preserving, hence the Veech group is a subgroup of $\SL(2,\mathbb{R})$.
For a finite translation surface, the Veech group is even a discrete subgroup of $\SL(2,\mathbb{R})$ but this fact relies heavily on the property that all singularities are cone angle singularities.

In these notes, we will subscribe to Chamanara's point of view that the Veech group is in fact the image of $\GL^+(X, \mathcal{A})$ in $\PSL(2,\mathbb{R})$.

\section{Definition of Chamanara surface}

In \cite{chamanara_04}, Chamanara describes in detail a family of translation surfaces with a parameter $\alpha \in (0,1)$. We will concentrate on the case of $\alpha = \frac{1}{2}$.
In general, Chamanara considers a square that has edges of length $\sum_{i=1}^\infty \alpha^{n}$ for an $\alpha \in (0,1)$ and the edges are divided into segments where the $n$th segment has length $\alpha^n$.

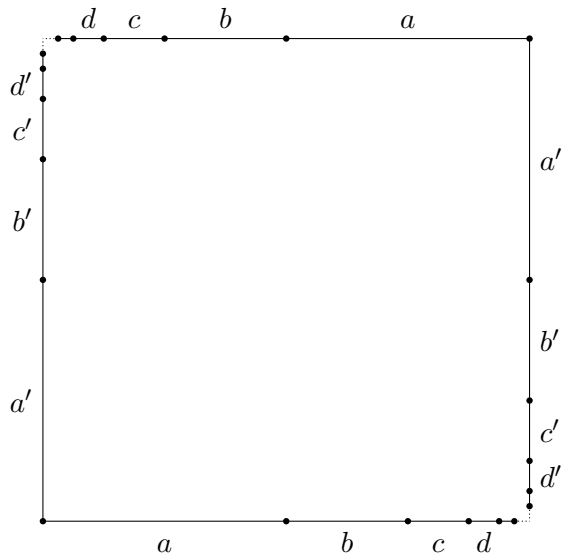
\begin{figure}[bthp]
 \begin{center}
 \begin{tikzpicture}[x=1cm,y=1cm, scale=0.8]
  \draw (0,0) -- node[below, text height=\heightof{$b$}] {$a$}
        (4,0) -- node[below, text height=\heightof{$b$}] {$b$}
        (6,0) -- node[below, text height=\heightof{$b$}] {$c$}
        (7,0) -- node[below, text height=\heightof{$b$}] {$d$}
        (7.5,0) -- (7.75,0);
  \draw[densely dotted] (7.75,0) -- (8,0);

  \fill (0,0) circle (1.5pt);
  \fill (8,8) circle (1.5pt);
  
  \fill (4,0) circle (1.5pt);
  \fill (6,0) circle (1.5pt);
  \fill (7,0) circle (1.5pt);
  \fill (7.5,0) circle (1.5pt);
  \fill (7.75,0) circle (1.5pt);

  \draw (8,8) -- node[above] {$a$}
        (4,8) -- node[above] {$b$}
        (2,8) -- node[above] {$c$}
        (1,8) -- node[above] {$d$}
        (0.5,8) -- (0.25,8);
  \draw[densely dotted] (0.25,8) -- (0,8);

  \fill (4,8) circle (1.5pt);
  \fill (2,8) circle (1.5pt);
  \fill (1,8) circle (1.5pt);
  \fill (0.5,8) circle (1.5pt);
  \fill (0.25,8) circle (1.5pt);

  \draw (0,0) -- node[left] {$a'$}
        (0,4) -- node[left] {$b'$}
        (0,6) -- node[left] {$c'$}
        (0,7) -- node[left] {$d'$}
        (0,7.5) -- (0,7.75);
  \draw[densely dotted] (0,7.75) -- (0,8);

  \fill (0,4) circle (1.5pt);
  \fill (0,6) circle (1.5pt);
  \fill (0,7) circle (1.5pt);
  \fill (0,7.5) circle (1.5pt);
  \fill (0,7.75) circle (1.5pt);

  \draw (8,8) -- node[right] {$a'$}
        (8,4) -- node[right] {$b'$}
        (8,2) -- node[right] {$c'$}
        (8,1) -- node[right] {$d'$}
        (8,0.5) -- (8,0.25);
  \draw[densely dotted] (8,0.25) -- (8,0);

  \fill (8,4) circle (1.5pt);
  \fill (8,2) circle (1.5pt);
  \fill (8,1) circle (1.5pt);
  \fill (8,0.5) circle (1.5pt);
  \fill (8,0.25) circle (1.5pt);
 \end{tikzpicture}
 \end{center}
 \caption{For the Chamanara surface, we identify segments that are parallel and have the same length.}
 \label{fig_chamanara}
\end{figure}

In our case, we consider a square with edges of length $1$. We divide the top edge into two halves and the bottom edge into two halves. Then we glue the right half of the top to the left half of the bottom. The remaining halves are divided again and the right part of the top is glued to the left part of the bottom, and so on (see \autoref{fig_chamanara}).
We do the same with the left and the right edge, always identifying the upper part of the right edge with the lower part of the left edge.

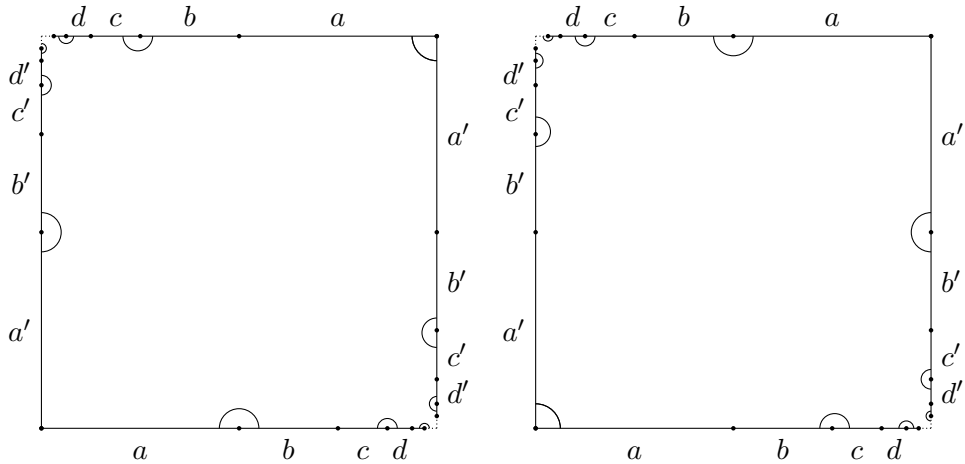
\begin{figure}
 \begin{center}
 \begin{tikzpicture}[scale=0.65, x=1cm,y=1cm]
  \newcommand\chamanaraseitea{
   \draw (0,0) -- (7.75,0);
   \draw[densely dotted] (7.75,0) -- (8,0);
   \fill (0,0) circle (1.3pt);
   \draw (3.6,0) arc (180:0:0.4);
   \fill (4,0) circle (1.3pt);
   \fill (6,0) circle (1.3pt);
   \draw (6.8,0) arc (180:0:0.2);
   \fill (7,0) circle (1.3pt);
   \fill (7.5,0) circle (1.3pt);
   \draw (7.65,0) arc (180:0:0.1);
   \fill (7.75,0) circle (1.3pt);
  }
  \newcommand\chamanaraseiteb{
   \draw (0,0) -- (7.75,0);
   \draw[densely dotted] (7.75,0) -- (8,0);
   \draw (0,0.5) arc (90:0:0.5);
   \fill (0,0) circle (1.3pt);
   \fill (4,0) circle (1.3pt);
   \draw (5.75,0) arc (180:0:0.3);
   \fill (6,0) circle (1.3pt);
   \fill (7,0) circle (1.3pt);
   \draw (7.35,0) arc (180:0:0.15);
   \fill (7.5,0) circle (1.3pt);
   \fill (7.75,0) circle (1.3pt);
  }
  \chamanaraseitea
  \path (0,0) -- node[below, text height=\heightof{$b$}] {$a$} (4,0) -- node[below, text height=\heightof{$b$}] {$b$} (6,0) -- node[below, text height=\heightof{$b$}] {$c$} (7,0) -- node[below, text height=\heightof{$b$}] {$d$} (7.5,0);

  \begin{scope}[xscale=-1,yscale=-1,xshift=-8cm,yshift=-8cm]
   \chamanaraseiteb
   \path (0,0) -- node[above] {$a$} (4,0) -- node[above] {$b$} (6,0) -- node[above] {$c$} (7,0) -- node[above] {$d$} (7.5,0);
  \end{scope}

  \begin{scope}[xscale=-1,rotate=90]
   \chamanaraseitea
   \path (0,0) -- node[left] {$a'$} (4,0) -- node[left] {$b'$} (6,0) -- node[left] {$c'$} (7,0) -- node[left] {$d'$} (7.5,0);
  \end{scope}

  \begin{scope}[xscale=-1,rotate=-90,xshift=-8cm,yshift=-8cm]
   \chamanaraseiteb
   \path (0,0) -- node[right] {$a'$} (4,0) -- node[right] {$b'$} (6,0) -- node[right] {$c'$} (7,0) -- node[right] {$d'$} (7.5,0);
  \end{scope}

  \begin{scope}[xshift=10cm]
   \chamanaraseiteb
   \path (0,0) -- node[below, text height=\heightof{$b$}] {$a$} (4,0) -- node[below, text height=\heightof{$b$}] {$b$} (6,0) -- node[below, text height=\heightof{$b$}] {$c$} (7,0) -- node[below, text height=\heightof{$b$}] {$d$} (7.5,0);

   \begin{scope}[xscale=-1,yscale=-1,xshift=-8cm,yshift=-8cm]
    \chamanaraseitea
    \path (0,0) -- node[above] {$a$} (4,0) -- node[above] {$b$} (6,0) -- node[above] {$c$} (7,0) -- node[above] {$d$} (7.5,0);
   \end{scope}

   \begin{scope}[xscale=-1,rotate=90]
    \chamanaraseiteb
    \path (0,0) -- node[left] {$a'$} (4,0) -- node[left] {$b'$} (6,0) -- node[left] {$c'$} (7,0) -- node[left] {$d'$} (7.5,0);
   \end{scope}

   \begin{scope}[xscale=-1,rotate=-90,xshift=-8cm,yshift=-8cm]
    \chamanaraseitea
    \path (0,0) -- node[right] {$a'$} (4,0) -- node[right] {$b'$} (6,0) -- node[right] {$c'$} (7,0) -- node[right] {$d'$} (7.5,0);
   \end{scope}
  \end{scope}
 \end{tikzpicture}
 \end{center}
 \label{fig_chamanara_singularities}
 \caption{The arcs indicate which cutting points are identified by the gluings in the Chamanara surface.}
\end{figure}

When we exclude not only the corners of the square but also the points on the edges where we divided the segments into halves, we obtain a translation structure on the resulting surface. The corners and the cutting points lead to points in the metric completion and so they are specifying singularities. By following the gluings we find that every second cutting point is identified as sketched in \autoref{fig_chamanara_singularities}.
However, the distance of the two indicated points is not bounded away from $0$ in the metric completion, hence it is~$0$. This means that all cutting points are identified to one point in the metric completion.
The same argument holds true for the corners and we receive that we have exactly one singularity $\sigma$.

\section{Cylinder decompositions and parabolic elements} \label{sec_cylinder_decompositions}

For a translation surface $(X,\mathcal{A})$ with at least one singularity, a \emph{cylinder} in $(X,\mathcal{A})$ of \emph{circumference} $w>0$ and \emph{height} $h>0$ is an open subset of $X$ which is isometric to a Euclidean cylinder $\mathbb{R}/w\mathbb{Z} \times (0,h)$.
The \emph{modulus} of a cylinder is the ratio of circumference and height, i.e.\ it is equal to $\frac{w}{h}$.

If a cylinder can be extended to a maximal cylinder then the maximal cylinder is bounded by geodesic segments in $X$ with endpoints in singularities. Segments of this kind are called \emph{saddle connections}. The \emph{direction} of the cylinder is the direction of the saddle connections.

A \emph{cylinder decomposition} of $(X,\mathcal{A})$ is a collection of maximal cylinders in $(X,\mathcal{A})$ so that the closures of the cylinders in $X$ cover $X$ and so that each two cylinders are disjoint.

We will use the following well-known result to determine elements of the Veech group. It was first proven in \cite{veech_89} for finite translation surfaces. However, the proof of statement~(i) literally works for the general case.

\begin{prop}[Cylinder decompositions and parabolic elements {[Veech]}] \label{prop_cylinder_decomposition_parabolic_element}
\leavevmode \vspace{-\baselineskip}
 \begin{enumerate}
  \item \label{item_cylinder_to_parabolic} Let $(X,\mathcal{A})$ be a translation surface and $\left\{z_n\right\}$ a cylinder decomposition so that the cylinder $z_n$ has height $h_n$ and circumference $w_n$. If the inverse moduli $\frac{h_n}{w_n}$ are commensurable, i.e.\ if there exists an $m\in \mathbb{R}$ so that each inverse modulus is an integer multiple of $m$, then the Veech group contains a parabolic element conjugated to the matrix
  $\begin{pmatrix}
   1 & \frac{1}{m} \\
   0 & 1  
  \end{pmatrix}$.

  \item \label{item_parabolic_to_cylinder} Let $(X,\mathcal{A})$ be a finite translation surface such that the Veech group contains a parabolic element. Then there exists a cylinder decomposition of $(X,\mathcal{A})$ in the eigen direction of the parabolic element.
 \end{enumerate}
\end{prop}

\begin{figure}[bt]
 \begin{center}
 \begin{tikzpicture}[x=1cm,y=1cm,scale=0.8]
  \newcommand\chamanaraseite{
   \draw (0,0) -- (7.75,0);
   \draw[densely dotted] (7.75,0) -- (8,0);
   \fill (0,0) circle (1pt);
   \fill (4,0) circle (1pt);
   \fill (6,0) circle (1pt);
   \fill (7,0) circle (1pt);
   \fill (7.5,0) circle (1pt);
   \fill (7.75,0) circle (1pt);
  }

  \draw[pattern color=gray!60, pattern=dots] (0,0) -- (2,0) -- (4,8) -- (2,8) -- (0,0);
  \draw[pattern color=gray!60, pattern=dots] (2,0) -- (4,0) -- (6,8) -- (4,8) -- (2,0);
  \draw[pattern color=gray!60, pattern=dots] (4,0) -- (6,0) -- (8,8) -- (6,8) -- (4,0);

  \draw[pattern color=gray!60, pattern=bricks] (0,0) -- (2,8) -- (1,8) -- (0,4) -- (0,0);
  \draw[pattern color=gray!60, pattern=bricks] (6,0) -- (7,0) -- (8,4) -- (8,8) -- (6,0);

  \draw[pattern color=gray!60, pattern=north west lines] (0,4) -- (1,8) -- (0.5,8) -- (0,6) -- (0,4);
  \draw[pattern color=gray!60, pattern=north west lines] (7,0) -- (7.5,0) -- (8,2) -- (8,4) -- (7,0);

  \draw[pattern color=gray!60, pattern=grid] (0,6) -- (0.5,8) -- (0.25,8) -- (0,7) -- (0,6);
  \draw[pattern color=gray!60, pattern=grid] (7.5,0) -- (7.75,0) -- (8,1) -- (8,2) -- (7.5,0);

  \chamanaraseite
  \path (0,0) -- node[below=0.11cm] {$a$} (4,0) -- node[below] {$b$} (6,0) -- node[below=0.11cm] {$c$} (7,0) -- node[below] {$d$} (7.5,0);

  \begin{scope}[xscale=-1,yscale=-1,xshift=-8cm,yshift=-8cm]
   \chamanaraseite
   \path (0,0) -- node[above] {$a$} (4,0) -- node[above] {$b$} (6,0) -- node[above] {$c$} (7,0) -- node[above] {$d$} (7.5,0);
  \end{scope}

  \begin{scope}[xscale=-1,rotate=90]
   \chamanaraseite
   \path (0,0) -- node[left] {$a'$} (4,0) -- node[left] {$b'$} (6,0) -- node[left] {$c'$} (7,0) -- node[left] {$d'$} (7.5,0);
  \end{scope}

  \begin{scope}[xscale=-1,rotate=-90,xshift=-8cm,yshift=-8cm]
   \chamanaraseite
   \path (0,0) -- node[right] {$a'$} (4,0) -- node[right] {$b'$} (6,0) -- node[right] {$c'$} (7,0) -- node[right] {$d'$} (7.5,0);
  \end{scope}
 \end{tikzpicture}
 \label{fig_cylinder_decomposition_chamanara_slope_4}
 \caption{Chamanara surface with a cylinder decomposition of slope $4$.}
 \end{center}
\end{figure}
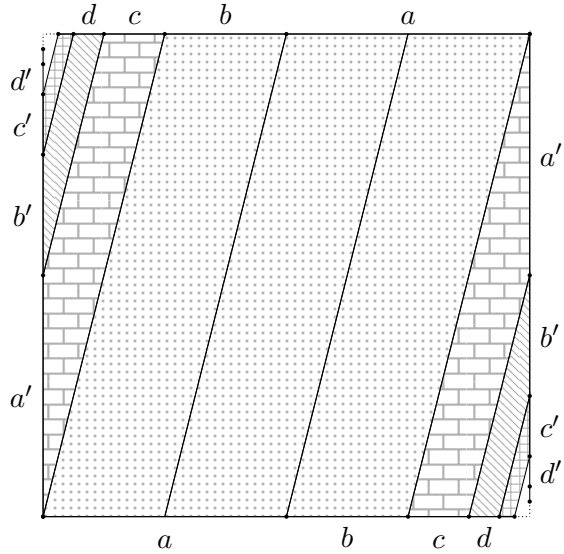

\autoref{fig_cylinder_decomposition_chamanara_slope_4} shows a cylinder decomposition of slope $4$ on the Chamanara surface. There are infinitely many cylinders where all except one can be seen as built by two trapezoids. As all trapezoids are similar, the moduli of these cylinders coincide. In fact, this modulus can be calculated to be $\frac{4}{51}$ by trigonometry.
Similar calculations show that the modulus of the middle cylinder is $\frac{4}{51}$, too. Then \autoref{prop_cylinder_decomposition_parabolic_element} implies that there exists a parabolic element in the Veech group which maps every saddle connection of slope $4$ to itself and acts as a Dehn twist on every cylinder in the cylinder decomposition.

\begin{figure}
 \begin{center}
 \begin{tikzpicture}[x=1cm,y=1cm,scale=0.8]
  \newcommand\chamanaraseite{
   \draw (0,0) -- (7.75,0);
   \draw[densely dotted] (7.75,0) -- (8,0);
   \fill (0,0) circle (1.5pt);
   \fill (4,0) circle (1.5pt);
   \fill (6,0) circle (1.5pt);
   \fill (7,0) circle (1.5pt);
   \fill (7.5,0) circle (1.5pt);
   \fill (7.75,0) circle (1.5pt);
  }

  \chamanaraseite
  \path (0,0) -- node[below=0.11cm] {$a$} (4,0) -- node[below] {$b$} (6,0) -- node[below=0.11cm] {$c$} (7,0) -- node[below] {$d$} (7.5,0);

  \begin{scope}[xscale=-1,yscale=-1,xshift=-8cm,yshift=-8cm]
   \chamanaraseite
   \path (0,0) -- node[above] {$a$} (4,0) -- node[above] {$b$} (6,0) -- node[above] {$c$} (7,0) -- node[above] {$d$} (7.5,0);
  \end{scope}

  \begin{scope}[xscale=-1,rotate=90]
   \chamanaraseite
   \path (0,0) -- node[left] {$a'$} (4,0) -- node[left] {$b'$} (6,0) -- node[left] {$c'$} (7,0) -- node[left] {$d'$} (7.5,0);
  \end{scope}

  \begin{scope}[xscale=-1,rotate=-90,xshift=-8cm,yshift=-8cm]
   \chamanaraseite
   \path (0,0) -- node[right] {$a'$} (4,0) -- node[right] {$b'$} (6,0) -- node[right] {$c'$} (7,0) -- node[right] {$d'$} (7.5,0);
  \end{scope}

  \draw[ultra thin] (0,7.75) -- (0.5,8);
  \draw[ultra thin] (0,7.5) -- (1,8);
  \draw[ultra thin] (0,7) -- (2,8);
  \draw[ultra thin] (0,6) -- (4,8);
  \draw[] (0,4) -- node[below, pos=0.3] {$\frac{2}{3} w$} (8,8);
  \draw[ultra thin] (0,0) -- (8,4);
  \draw[] (4,0) -- node[above] {$\frac{1}{3} w$} (8,2);
  \draw[ultra thin] (6,0) -- (8,1);
  \draw[ultra thin] (7,0) -- (8,0.5);
  \draw[ultra thin] (7.5,0) -- (8,0.25);

  \draw (3,7.5) -- node[right] {$h$} (3.8,5.9);
 \end{tikzpicture}
 \label{fig_chamanara_cylinder_decomposition_with_values}
 \caption{Chamanara surface with saddle connections of slope $\frac{1}{2}$: $h$ is the height and $w$ is the circumference of the largest cylinder that is built by two trapezoids.}
 \end{center}
\end{figure}
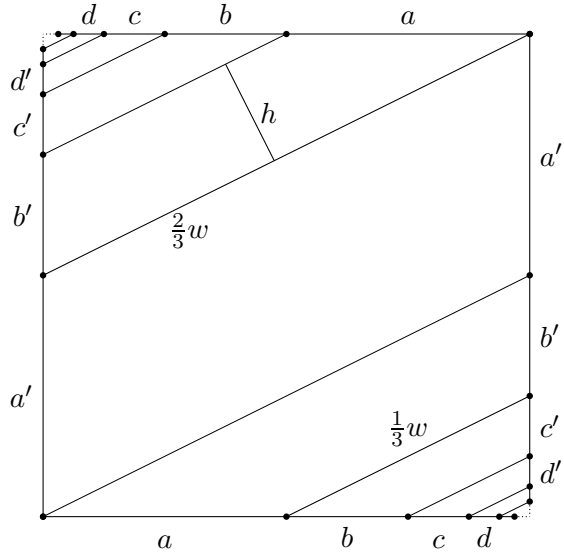

For every slope $2^n$ with $n\in \mathbb{Z}$ there exists a cylinder decomposition like the previous one (cf.\ \cite[Proposition 11]{chamanara_04}):
We have that every maximal geodesic of slope $2^n$ starting in a singularity is a saddle connection and bounds a cylinder of the same slope.
As in the previous case, we have infinitely many cylinders built by two trapezoids and finitely many cylinders built by parallelograms like the middle cylinder in the previous case.
Again, the moduli of the cylinders of the first type are all the same. 
Let $h$ be the height of the largest cylinder of the first type and $w$ its circumference (see \autoref{fig_chamanara_cylinder_decomposition_with_values}). Then we have for a cylinder of the second type that its height is an integer multiple of $h$ and its circumference is an integer multiple of $\frac{2}{3} w$. Therefore the inverse modulus of this cylinder is a rational multiple of $\frac{3}{2} \cdot \frac{h}{w}$.
For example, for slope $2$ and slope $\frac{1}{2}$, the inverse modulus of the middle cylinder is three times the other inverse moduli.
This yields the existence of a number $m$ such that for every cylinder $z_i$ there exists an integer $k_i$ such that the inverse modulus of $z_i$ is $k_i \cdot m$. By \autoref{prop_cylinder_decomposition_parabolic_element}, there exists a parabolic element in the Veech group such that the corresponding affine map twists the $i$th cylinder $k_i$ times.

\bigskip

We proceed by explicitly calculating two parabolic elements in the Veech group. These correspond to the Dehn twists in the cylinder decompositions of slope $1$ and slope $2$.

\begin{figure}[p]
 \begin{center}
 \begin{tikzpicture}[x=1cm,y=1cm,scale=0.8]
  \newcommand\chamanaraseite{
   \draw (0,0) -- (7.75,0);
   \draw[densely dotted] (7.75,0) -- (8,0);
   \fill (0,0) circle (1pt);
   \fill (4,0) circle (1pt);
   \fill (6,0) circle (1pt);
   \fill (7,0) circle (1pt);
   \fill (7.5,0) circle (1pt);
   \fill (7.75,0) circle (1pt);
  }

  \draw[pattern color=gray!60, pattern=dots] (0,0) -- (8,8) -- (4,8) -- (0,4) -- (0,0);
  \draw[pattern color=gray!60, pattern=dots] (0,0) -- (4,0) -- (8,4) -- (8,8) -- (0,0);

  \draw[pattern color=gray!60, pattern=bricks] (0,4) -- (4,8) -- (2,8) -- (0,6) -- (0,4);
  \draw[pattern color=gray!60, pattern=bricks] (4,0) -- (6,0) -- (8,2) -- (8,4) -- (4,0);

  \draw[pattern color=gray!60, pattern=north west lines] (0,6) -- (2,8) -- (1,8) -- (0,7) -- (0,6);
  \draw[pattern color=gray!60, pattern=north west lines] (6,0) -- (7,0) -- (8,1) -- (8,2) -- (6,0);

  \draw[pattern color=gray!60, pattern=grid] (0,7) -- (1,8) -- (0.5,8) -- (0,7.5) -- (0,7);
  \draw[pattern color=gray!60, pattern=grid] (7,0) -- (7.5,0) -- (8,0.5) -- (8,1) -- (7,0);
  
  \chamanaraseite
  \path (0,0) -- node[below=0.11cm] {$a$} (4,0) -- node[below] {$b$} (6,0) -- node[below=0.11cm] {$c$} (7,0) -- node[below] {$d$} (7.5,0);

  \begin{scope}[xscale=-1,yscale=-1,xshift=-8cm,yshift=-8cm]
   \chamanaraseite
   \path (0,0) -- node[above] {$a$} (4,0) -- node[above] {$b$} (6,0) -- node[above] {$c$} (7,0) -- node[above] {$d$} (7.5,0);
  \end{scope}

  \begin{scope}[xscale=-1,rotate=90]
   \chamanaraseite
   \path (0,0) -- node[left] {$a'$} (4,0) -- node[left] {$b'$} (6,0) -- node[left] {$c'$} (7,0) -- node[left] {$d'$} (7.5,0);
  \end{scope}

  \begin{scope}[xscale=-1,rotate=-90,xshift=-8cm,yshift=-8cm]
   \chamanaraseite
   \path (0,0) -- node[right] {$a'$} (4,0) -- node[right] {$b'$} (6,0) -- node[right] {$c'$} (7,0) -- node[right] {$d'$} (7.5,0);
  \end{scope}
 \end{tikzpicture}
 \label{fig_cylinder_decomposition_chamanara_slope_1}
 \caption{Chamanara surface with a cylinder decomposition of slope $1$.}
 \end{center}
\end{figure}
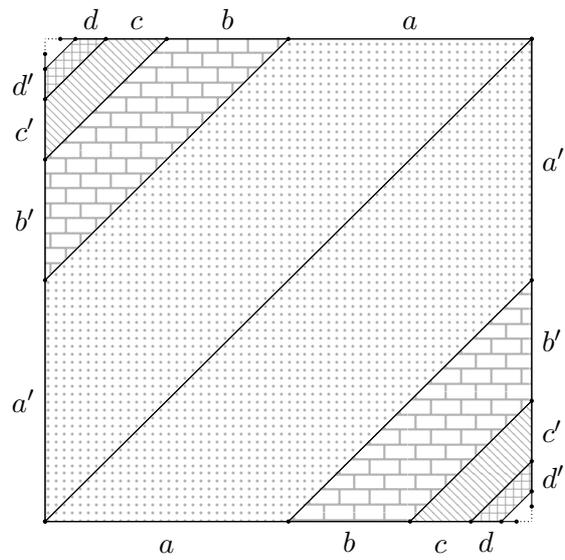

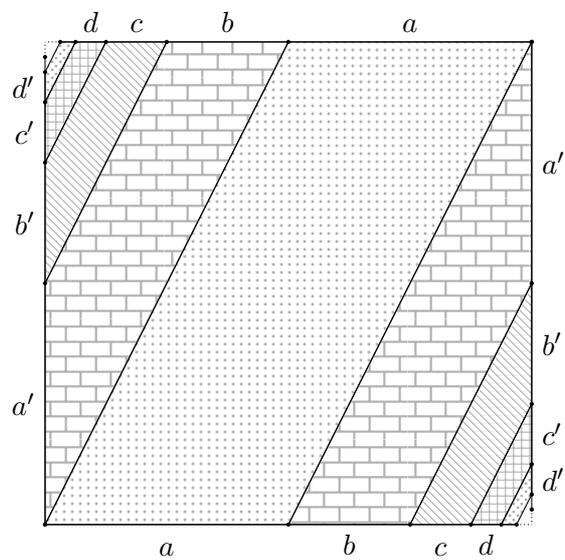
\begin{figure}[p]
 \begin{center}
 \begin{tikzpicture}[x=1cm,y=1cm,scale=0.8]
  \newcommand\chamanaraseite{
   \draw (0,0) -- (7.75,0);
   \draw[densely dotted] (7.75,0) -- (8,0);
   \fill (0,0) circle (1pt);
   \fill (4,0) circle (1pt);
   \fill (6,0) circle (1pt);
   \fill (7,0) circle (1pt);
   \fill (7.5,0) circle (1pt);
   \fill (7.75,0) circle (1pt);
  }

  \draw[pattern color=gray!60, pattern=dots] (0,0) -- (4,0) -- (8,8) -- (4,8) -- (0,0);

  \draw[pattern color=gray!60, pattern=bricks] (0,0) -- (4,8) -- (2,8) -- (0,4) -- (0,0);
  \draw[pattern color=gray!60, pattern=bricks] (4,0) -- (6,0) -- (8,4) -- (8,8) -- (4,0);

  \draw[pattern color=gray!60, pattern=north west lines] (0,4) -- (2,8) -- (1,8) -- (0,6) -- (0,4);
  \draw[pattern color=gray!60, pattern=north west lines] (6,0) -- (7,0) -- (8,2) -- (8,4) -- (6,0);

  \draw[pattern color=gray!60, pattern=grid] (0,6) -- (1,8) -- (0.5,8) -- (0,7) -- (0,6);
  \draw[pattern color=gray!60, pattern=grid] (7,0) -- (7.5,0) -- (8,1) -- (8,2) -- (7,0);
  
  \draw[pattern color=gray!60, pattern=crosshatch dots] (0,7) -- (0.5,8) -- (0.25,8) -- (0,7.5) -- (0,7);
  \draw[pattern color=gray!60, pattern=crosshatch dots] (7.5,0) -- (7.75,0) -- (8,0.5) -- (8,1) -- (7.5,0);

  \chamanaraseite
  \path (0,0) -- node[below=0.11cm] {$a$} (4,0) -- node[below] {$b$} (6,0) -- node[below=0.11cm] {$c$} (7,0) -- node[below] {$d$} (7.5,0);

  \begin{scope}[xscale=-1,yscale=-1,xshift=-8cm,yshift=-8cm]
   \chamanaraseite
   \path (0,0) -- node[above] {$a$} (4,0) -- node[above] {$b$} (6,0) -- node[above] {$c$} (7,0) -- node[above] {$d$} (7.5,0);
  \end{scope}

  \begin{scope}[xscale=-1,rotate=90]
   \chamanaraseite
   \path (0,0) -- node[left] {$a'$} (4,0) -- node[left] {$b'$} (6,0) -- node[left] {$c'$} (7,0) -- node[left] {$d'$} (7.5,0);
  \end{scope}

  \begin{scope}[xscale=-1,rotate=-90,xshift=-8cm,yshift=-8cm]
   \chamanaraseite
   \path (0,0) -- node[right] {$a'$} (4,0) -- node[right] {$b'$} (6,0) -- node[right] {$c'$} (7,0) -- node[right] {$d'$} (7.5,0);
  \end{scope}
 \end{tikzpicture}
 \label{fig_cylinder_decomposition_chamanara_slope_2}
 \caption{Chamanara surface with a cylinder decomposition of slope $2$.}
 \end{center}
\end{figure}

For slope $1$, all cylinders are built by two trapezoids (see \autoref{fig_cylinder_decomposition_chamanara_slope_1}). For the largest cylinder, the height is $\frac{\sqrt{2}}{4}$ and the circumference is $\sqrt{2} + \frac{\sqrt{2}}{2} = 3 \cdot \frac{\sqrt{2}}{2}$. Therefore, the modulus of all cylinders with slope $1$ is $\frac{3}{2} \cdot 4 = 6$.

For slope $2$, all except one cylinder are built by trapezoids (see \autoref{fig_cylinder_decomposition_chamanara_slope_2}). For the largest of these cylinders, the height is $\frac{1}{2 \sqrt{5}}$ and the circumference is $\frac{\sqrt{5}}{2} + \frac{\sqrt{5}}{4} = 3 \cdot \frac{\sqrt{5}}{4}$. Therefore, the modulus of all cylinders with slope $2$ of this type is $\frac{15}{2}$.
The height of the middle cylinder is $\frac{1}{\sqrt{5}}$ and the circumference is $\frac{\sqrt{5}}{2}$, hence the modulus is $\frac{5}{2}$.

For the sake of convenience in calculations, we will now rotate the Chamanara surface by $-\frac{\pi}{4}$ so that the former cylinder decomposition of slope $1$ is now horizontal.

Thus, by \autoref{prop_cylinder_decomposition_parabolic_element}, the element
$\begin{pmatrix}
  1 & 6 \\
  0 & 1
 \end{pmatrix}$,
which acts as a Dehn twist on the horizontal cylinders, is contained in the Veech group.
Note that for every parabolic element in the Veech group, for which the eigen direction is the horizontal direction, all horizontal saddle connections have to be fixed pointwise or are inverted in the center of the square. In particular, such a parabolic element has to act as a (multiple) Dehn twist on the horizontal cylinders, so up to rotation by
$\begin{pmatrix}
 -1 & 0 \\
  0 & -1
\end{pmatrix}$ it is of the form
$\begin{pmatrix}
 1 & 6k \\
 0 & 1
\end{pmatrix}$ for some $k\in \mathbb{Z} \setminus \{0\}$.

The second of the discussed cylinder decompositions has an angle $\alpha$ with the horizontal axis for which it holds $\sin \alpha = \frac{\sqrt{2}}{4} \cdot \frac{2}{\sqrt{5}} = \frac{1}{\sqrt{10}}$. Therefore, we know that the following parabolic matrix is contained in the Veech group.
\begin{align*}
 &
 \begin{pmatrix}
  \cos \alpha & - \sin \alpha \\
  \sin \alpha & \cos \alpha \\
 \end{pmatrix}
 \cdot
 \begin{pmatrix}
  1 & \frac{15}{2} \\
  0 & 1
 \end{pmatrix}
 \cdot
 \begin{pmatrix}
  \cos -\alpha & - \sin -\alpha \\
  \sin -\alpha & \cos -\alpha \\
 \end{pmatrix}
 \\
 = &
 \begin{pmatrix}
  \frac{3}{\sqrt{10}} & - \frac{1}{\sqrt{10}} \\
  \frac{1}{\sqrt{10}} & \frac{3}{\sqrt{10}} \\
 \end{pmatrix}
 \cdot
 \begin{pmatrix}
  1 & \frac{15}{2} \\
  0 & 1
 \end{pmatrix}
 \cdot
 \begin{pmatrix}
  \frac{3}{\sqrt{10}} & \frac{1}{\sqrt{10}} \\
  - \frac{1}{\sqrt{10}} & \frac{3}{\sqrt{10}} \\
 \end{pmatrix}
 \\
 = &
 \frac{1}{10} \cdot
 \begin{pmatrix}
  3 & \frac{43}{2} \\
  1 & \frac{21}{2}
 \end{pmatrix}
 \cdot
 \begin{pmatrix}
  3 & 1 \\
  -1 & 3 \\
 \end{pmatrix}
 = 
 \frac{1}{4} \cdot
 \begin{pmatrix}
  -5 & 27 \\
  -3 & 13
 \end{pmatrix}
\end{align*}

This parabolic element acts as a triple Dehn twist on the middle cylinder and as a single Dehn twist on the other cylinders. So, with the same argument as before, a parabolic element with the same eigen direction has to be a power of $\frac{1}{4} \cdot
 \begin{pmatrix}
  -5 & 27 \\
  -3 & 13
 \end{pmatrix}$ -- up to rotation by
$\begin{pmatrix}
 -1 & 0 \\
  0 & -1
\end{pmatrix}$.

We summarize this in the following lemma.

\begin{lem}\label{lem_P1_and_P2_are_maximal}
 A parabolic element in the Veech group $\Gamma$ of the Chamanara surface
 \begin{enumerate}
  \item with horizontal eigen direction is a power of $P_1 \coloneqq
  \begin{pmatrix}
   1 & 6 \\
   0 & 1
  \end{pmatrix}$ and
  \item with eigen direction $\alpha$ is a power of $P_2 \coloneqq \frac{1}{4} \cdot
  \begin{pmatrix}
   -5 & 27 \\
   -3 & 13
  \end{pmatrix}$.
 \end{enumerate}
\end{lem}

To finish this section we give a large set of directions that can not occur as eigen directions of parabolic elements in the Veech group. Note that there can still exist cylinder decompositions in these directions.

\begin{lem}\label{lem_forbidden_eigen_directions}
 The eigen direction of a parabolic element in the Veech group $\Gamma$ of the Chamanara surface is contained in $[-\frac{\pi}{4}, \frac{\pi}{4}]$.
 
\begin{proof}
 Consider a direction $\theta$ which makes an angle of more than $\frac{\pi}{4}$ with the horizontal direction and suppose $P$ is a parabolic element with eigen direction $\theta$. There exist infinitely many geodesic segments in direction $\theta$ that start in the singularity and have length $1$ (see \autoref{fig_non_parabolic_directions}).
 Then $P$ has to map the set of these geodesic segments to itself. Note that the geodesic segment that starts in the upper corner is distinguished from all others as it is the limit of the geodesic segments in this set. Hence, $P$ has to fix it pointwise, up to a rotation by 
 $\begin{pmatrix}
  -1 & 0 \\
   0 & -1
 \end{pmatrix}$.
 As $P$ is locally acting as a shear on the Chamanara surface, this implies that also the geodesic segments close to the limit segment have to be fixed pointwise.
 So $P$ fixes an open subset of the Chamanara surface. Therefore, $P$ is the identity and not a parabolic element.
\begin{figure}
 \begin{center}
 \begin{tikzpicture}[x=1cm,y=1cm,scale=0.8, rotate=-45]
  \newcommand\chamanaraseite{
   \draw (0,0) -- (7.75,0);
   \draw[densely dotted] (7.75,0) -- (8,0);
   \fill (0,0) circle (1pt);
   \fill (4,0) circle (1pt);
   \fill (6,0) circle (1pt);
   \fill (7,0) circle (1pt);
   \fill (7.5,0) circle (1pt);
   \fill (7.75,0) circle (1pt);
  }

  \chamanaraseite
  \path (0,0) -- node[below=0.11cm] {$a$} (4,0) -- node[below] {$b$} (6,0) -- node[below=0.11cm] {$c$} (7,0) -- node[below] {$d$} (7.5,0);

  \begin{scope}[xscale=-1,yscale=-1,xshift=-8cm,yshift=-8cm]
   \chamanaraseite
   \path (0,0) -- node[above] {$a$} (4,0) -- node[above] {$b$} (6,0) -- node[above] {$c$} (7,0) -- node[above] {$d$} (7.5,0);
  \end{scope}

  \begin{scope}[xscale=-1,rotate=90]
   \chamanaraseite
   \path (0,0) -- node[left] {$a'$} (4,0) -- node[left] {$b'$} (6,0) -- node[left] {$c'$} (7,0) -- node[left] {$d'$} (7.5,0);
  \end{scope}

  \begin{scope}[xscale=-1,rotate=-90,xshift=-8cm,yshift=-8cm]
   \chamanaraseite
   \path (0,0) -- node[right] {$a'$} (4,0) -- node[right] {$b'$} (6,0) -- node[right] {$c'$} (7,0) -- node[right] {$d'$} (7.5,0);
  \end{scope}
  
  \draw[thick] (4,8) -- ++(-65:8cm);
  \draw[thick] (2,8) -- ++(-65:8cm);
  \draw[thick] (1,8) -- ++(-65:8cm);
  \draw[thick] (0.5,8) -- ++(-65:8cm);
  \draw[thick] (0.25,8) -- ++(-65:8cm);
  \draw[thick] (0.125,8) -- ++(-65:8cm);
  \draw[thick] (0.0625,8) -- ++(-65:8cm);
  \draw[thick] (0,8) -- ++(-65:8cm);
 \end{tikzpicture}
 \end{center}
 \caption{Geodesic segments in a direction which can not be the eigen direction of a parabolic element in the Veech group.}
 \label{fig_non_parabolic_directions}
\end{figure}
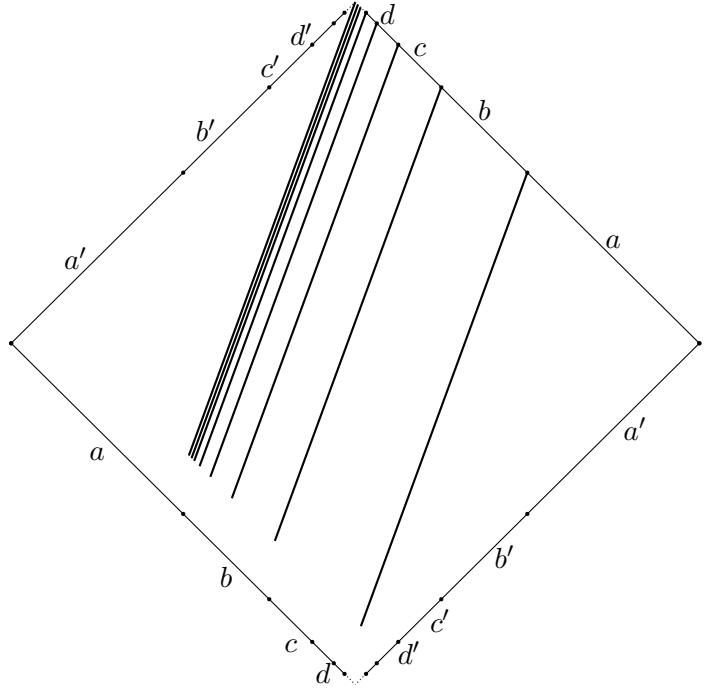
\end{proof}
\end{lem}

\section{A subgroup of the Veech group and its fundamental domain}

Consider the group $G$ which is the image of the group generated by the two parabolic elements
\begin{equation*}
 P_1 \coloneqq
 \begin{pmatrix}
  1 & 6 \\
  0 & 1
 \end{pmatrix},
 \quad
 P_2 \coloneqq
 \frac{1}{4} \cdot
 \begin{pmatrix}
  -5 & 27 \\
  -3 & 13
 \end{pmatrix}
\end{equation*}
in $\PSL(2,\mathbb{R})$. As we have seen in the last section, this is a subgroup of the Veech group. We will understand its fundamental domain in $\mathbb{H}$ in this section.

\bigskip

The elements of $\PSL(2,\mathbb{R})$ act as Möbius transformations on $\overline{\mathbb{H}} \coloneqq \mathbb{H} \cup\mathbb{R} \cup \infty$. Under this action, the fixed points of $P_1$ and $P_2$ are $\infty$ and $3$, respectively.

The product of the two generators is
\begin{equation*}
 H \coloneqq P_2 \cdot P_1 =
 \frac{1}{4} \cdot
 \begin{pmatrix}
  -5 & 27 \\
  -3 & 13
 \end{pmatrix}
 \cdot
 \begin{pmatrix}
  1 & 6 \\
  0 & 1
 \end{pmatrix}
 =
 -\frac{1}{4} \cdot
 \begin{pmatrix}
  5 & 3 \\
  3 & 5
 \end{pmatrix}
 ,
\end{equation*}
which is a hyperbolic element in $G$ with
$
H^{-1} = \frac{1}{4} \cdot
 \begin{pmatrix}
  -5 & 3 \\
  3 & -5
 \end{pmatrix}
$.
The element $H$ has fixed points $1$ and $-1$, in particular it fixes the hyperbolic geodesic through $-1$, $\mathrm{i}$, and $1$.

Of course, also $P_1$ and $H$ generate $G$. We will use these as generators as they make the calculations easier.

\begin{figure}[p]
 \begin{center}
  \begin{tikzpicture}[scale=1.85]
   \fill[gray!20!white] (-3,0) -- (3,0) -- (3,2.2) -- (-3,2.2) -- (-3,0);
   \draw (-3,0) -- (-3,2.2);
   \draw (3,0) -- (3,2.2);

   \draw[->, gray] (-3.3,0) -- (3.3,0) node[right] {};
   \draw[->, gray] (0,0) -- (0,2.3) node[above] {};		

   \draw (-3,0.05) -- (-3,-0.05) node[below] {$-3$};		
   \draw (3,0.05) -- (3,-0.05) node[below] {$3$};
   \draw (-1,0.05) -- (-1,-0.05) node[below] {$-1$};		
   \draw (1,0.05) -- (1,-0.05) node[below] {$1$};
   \draw (0,0.05) -- (0,-0.05) node[below] {$0$};
   \draw (-0.05,1) -- (0.05,1) node[right] {$\mathrm{i}$};				
   
   \draw (2,1.7) node {$F_1$};
  \end{tikzpicture}
 \end{center}
 \caption{Fundamental domain of $\left\langle P_1 \right\rangle$.}
 \label{fig_fundamental_domain_P_1}
\end{figure}
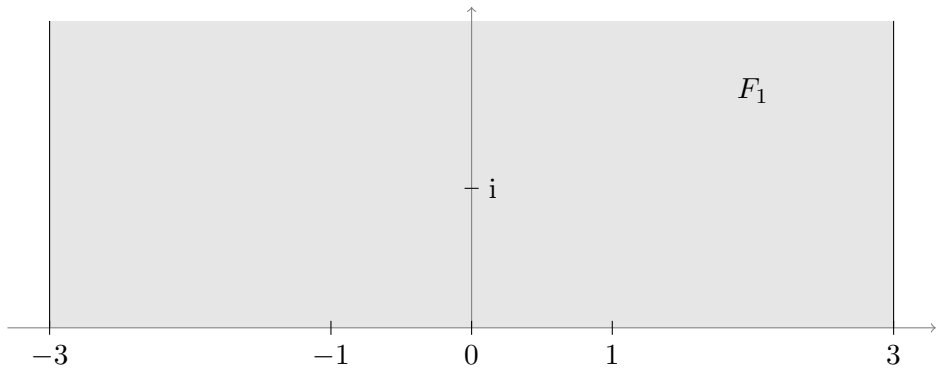

A fundamental domain of the subgroup $\left\langle P_1 \right\rangle$ is the strip $F_1 \coloneqq \{z \in \mathbb{H} : -3 < \Re(z) < 3\}$ as shown in \autoref{fig_fundamental_domain_P_1}.

We construct a fundamental domain of $\left\langle H \right\rangle$ by the Dirichlet method:
Let $M$ be the matrix
$\frac{\sqrt{2}}{2} \cdot
\begin{pmatrix}
 1 & -1 \\
 1 & 1
\end{pmatrix}$. This defines an elliptic element of order $4$ which fixes $\mathrm{i}$. We have
\begin{equation*}
 M^{-1}HM =
 - \frac{1}{4} \cdot \frac{1}{2}
 \begin{pmatrix}
  1 & 1 \\
  -1 & 1
 \end{pmatrix}
 \cdot
 \begin{pmatrix}
  5 & 3 \\
  3 & 5
 \end{pmatrix}
 \cdot
 \begin{pmatrix}
  1 & -1 \\
  1 & 1
 \end{pmatrix}
 =
 - \frac{1}{8} \cdot
 \begin{pmatrix}
  8 & 8 \\
  -2 & 2
 \end{pmatrix}
 \cdot
 \begin{pmatrix}
  1 & -1 \\
  1 & 1
 \end{pmatrix}
 =
 \begin{pmatrix}
  2 & 0 \\
  0 & \frac{1}{2}
 \end{pmatrix} .
\end{equation*}

\begin{figure}[p]
 \begin{center}
  \begin{tikzpicture}[scale=1.85]
   \fill[gray!20!white] (-3.2,0) -- (-3,0) arc (180:0:1.333) -- (0.333,0) arc (180:0:1.333) -- (3.2,0) -- (3.2,2.2) -- (-3.2,2.2) -- (-3.2,0);
   \draw (-3,0) arc (180:0:1.333);
   \draw (0.333,0) arc (180:0:1.333);

   \draw[->, gray] (-3.3,0) -- (3.3,0) node[right] {};
   \draw[->, gray] (0,0) -- (0,2.3) node[above] {};		

   \draw (-3,0.05) -- (-3,-0.05) node[below] {$-3$};		
   \draw (3,0.05) -- (3,-0.05) node[below] {$3$};
   \draw (-1,0.05) -- (-1,-0.05) node[below] {$-1$};		
   \draw (1,0.05) -- (1,-0.05) node[below] {$1$};
   \draw (-0.333,0.05) -- (-0.333,-0.05) node[below] {$-\frac{1}{3}$};		
   \draw (0.333,0.05) -- (0.333,-0.05) node[below] {$\frac{1}{3}$};
   \draw (0,0.05) -- (0,-0.05) node[below] {$0$};
   \draw (-0.05,1) -- (0.05,1) node[right] {$\mathrm{i}$};				
   
   \draw (2.5,1.7) node {$F_2$};
  \end{tikzpicture}
 \end{center}
 \caption{Fundamental domain of $\left\langle H \right\rangle$.}
 \label{fig_fundamental_domain_H}
\end{figure}
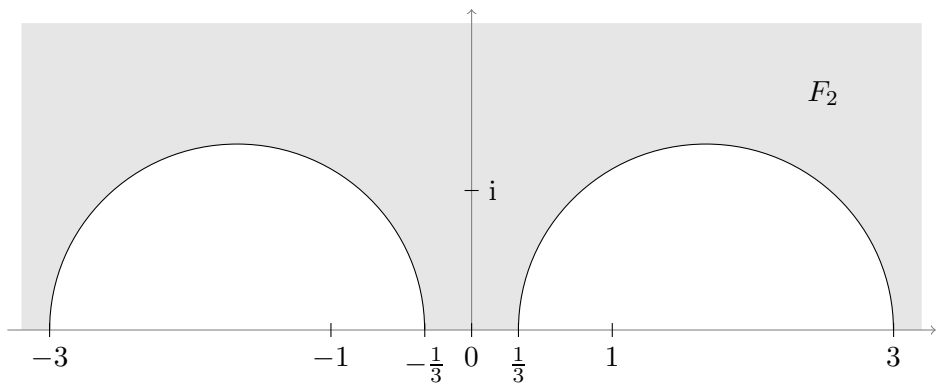

A fundamental domain for $\left\langle M^{-1}HM \right\rangle$ is the annulus $A = \{ z \in \mathbb{H} : \frac{1}{2} < |z| < 2 \}$. Now we have that the image of $A$ under $M$ is a fundamental domain for $\left\langle H \right\rangle$. This is even a Dirichlet fundamental domain with center $\mathrm{i}$ as $M$ fixes $\mathrm{i}$.
We have $M \bullet -2 = 3$, $M \bullet -\frac{1}{2} = -3$, $M \bullet \frac{1}{2} = -\frac{1}{3}$, and $M \bullet 2 = \frac{1}{3}$. Hence, we have that $F_2$ as in \autoref{fig_fundamental_domain_H} is a fundamental domain of $\left\langle H \right\rangle$.

\begin{figure}[p]
 \begin{center}
  \begin{tikzpicture}[scale=1.85]
   \fill[gray!20!white] (-3,0) arc (180:0:1.333) -- (0.333,0) arc (180:0:1.333) -- (3,2.2) -- (-3,2.2) -- (-3,0);
   \draw (-3,2.2) -- (-3,0) arc (180:0:1.333);
   \draw (0.333,0) arc (180:0:1.333) -- (3,2.2);

   \draw[->, gray] (-3.3,0) -- (3.3,0) node[right] {};
   \draw[->, gray] (0,0) -- (0,2.3) node[above] {};		

   \draw (-3,0.05) -- (-3,-0.05) node[below] {$-3$};		
   \draw (3,0.05) -- (3,-0.05) node[below] {$3$};
   \draw (-1,0.05) -- (-1,-0.05) node[below] {$-1$};		
   \draw (1,0.05) -- (1,-0.05) node[below] {$1$};
   \draw (-0.333,0.05) -- (-0.333,-0.05) node[below] {$-\frac{1}{3}$};		
   \draw (0.333,0.05) -- (0.333,-0.05) node[below] {$\frac{1}{3}$};
   \draw (0,0.05) -- (0,-0.05) node[below] {$0$};
   \draw (-0.05,1) -- (0.05,1) node[right] {$\mathrm{i}$};				
   
   \draw (2.5,1.7) node {$F$};
  \end{tikzpicture}
 \end{center}
 \caption{Candidate for the fundamental domain of $G$.}
 \label{fig_fundamental_domain_G}
\end{figure}
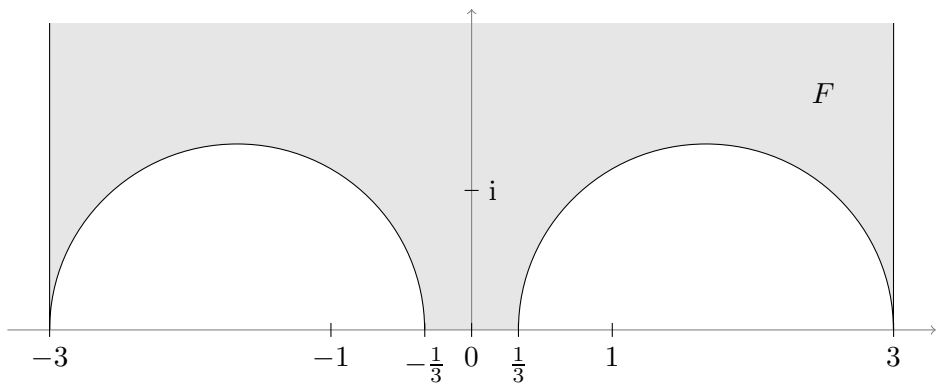

It follows that $F \coloneqq F_1 \cap F_2$ is of the form as sketched in \autoref{fig_fundamental_domain_G}.

It remains to show that $F$ is in fact a fundamental domain of $G$. To this end, we apply Poincaré's theorem (see \cite[Section 9.8]{beardon_95}). $F$ is an open convex hyperbolic polygon, $P_1$ and $H$ constitute a side pairing for $F$, and the cusps at the boundary of $F$, namely $\infty$, $3$, and $-3$, are fixed points of parabolic elements, namely $P_1$, $P_2$, and $P_1^{-1} P_2 P_1$.
Hence Poincaré's theorem tells us that the group generated by the side pairings, i.e.\ $G$, is discrete and that $F_G \coloneqq F$ is a fundamental domain for $G$.

Note that the cusps at $3$ and $-3$ are identified, so the group $G$ has two cusps and one hole.

\section{The Veech group}

Let us first note that the Veech group $\Gamma$ of the Chamanara surface is also discrete.
For this, we show that the limit set of $\Gamma$ has empty intersection with the open real segment from $-1$ to $1$.

First, we determine which eigen directions correspond to fixed points of parabolic elements in $(-1,1)$.
The parabolic element
$\frac{1}{2} \cdot
 \begin{pmatrix}
  1 & -1 \\
  1 & 3
 \end{pmatrix}
$
has eigen direction $-\frac{\pi}{4}$ and fixed point $-1 \in \partial \mathbb{H}$, whereas the parabolic element
$\frac{1}{2} \cdot
 \begin{pmatrix}
  1 & 1 \\
  -1 & 3
 \end{pmatrix}
$
has eigen direction $\frac{\pi}{4}$ and fixed point $1 \in \partial \mathbb{H}$.
Therefore, every parabolic element that has a fixed point in $(-1,1)$ has an eigen direction which makes an angle of more than $\frac{\pi}{4}$ with the horizontal direction.
We have already seen in \autoref{lem_forbidden_eigen_directions} that no parabolic element with such an eigen direction exists.

So, suppose that a point in $(-1,1)$ is contained in the limit set of $\Gamma$.
Since the fixed points of hyperbolic elements of $\Gamma$ are dense in the limit set (see the more general version for Kleinian groups in \cite[Proposition V.E.3.]{maskit_88}), there also exists a point $x \in (-1,1)$ which is the fixed point of a hyperbolic element.
In particular, $x$ is the attracting fixed point of a hyperbolic element $h \in \Gamma$. Then there exists an $n\in \mathbb{N}$ so that $h^n(\infty)$ is contained in $(-1,1)$. On the other hand, $h^n(\infty)$ is the fixed point of the parabolic element $h^n P_1 h^{-n}$. Again, as there are no parabolic elements in $\Gamma$ with a fixed point in $(-1,1)$, also $h$ can not be contained in $\Gamma$.

From this we can deduce the nature of the group $\Gamma$.

\begin{prop}
 The Veech group $\Gamma$ of the Chamanara surface is a non-elementary Fuchsian group of the second kind.

\begin{proof}
 As $\Gamma$ acts discontinuously on a neighborhood of $(-1,1)$ in $\mathbb{H}$, the group $\Gamma$ is a discrete subgroup of $\PSL(2,\mathbb{R})$, i.e.\ a Fuchsian group.
 
 As the limit set of $\Gamma$ is disjoint from $(-1,1)$, $\Gamma$ is of the second kind. On the other hand, the fixed points of $P_1$, $P_2$, and $P_1^{-1} P_2 P_1$ are different elements of the limit set. This implies that the limit set contains infinitely many elements and $\Gamma$ is not elementary.
\end{proof}
\end{prop}

We proceed with the observation that $G$ is in fact a normal subgroup of the Veech group $\Gamma$. For this, we prove that for any element $A\in \Gamma$ and for $i = 1,2$, the parabolic element $A^{-1} P_i A$ is contained in $G$.

Choose a fundamental domain $F_\Gamma$ of $\Gamma$ with $F_\Gamma \subseteq F_G$ and so that $F_G$ is tessellated by translates of $F_\Gamma$ under elements of $\Gamma$.
Let $A\in \Gamma$ and $i \in \{1,2\}$. Then $A^{-1} P_i A$ is parabolic and has a fixed point $c \in \mathbb{R} \cup \{\infty\}$. There exists an element $g \in G$ so that $g \bullet c$ is contained in the closure of $F_G$. It is the fixed point of the parabolic element $g A^{-1} P_i A g^{-1}$.
This gives us two possibilities for $g \bullet c$: either it has to be contained in the hole of $F_G$ or it is one of the three cusps of $F_G$. The first case is not possible as the points in the hole are contained in $(-1,1)$ and hence can not be fixed points of parabolic elements in the Veech group as we have noted before.
In the second case, $g \in G$ can be chosen so that $g \bullet c$ is $3$ or $\infty$. Therefore, $g A^{-1} P_i A g^{-1}$ has the same eigen direction as $P_1$ or $P_2$ and as shown in \autoref{lem_P1_and_P2_are_maximal} this means that $g A^{-1} P_i A g^{-1}$ is a power of $P_1$ or $P_2$. Hence, $g A^{-1} P_i A g^{-1}$ is contained in $G$ and so is $A^{-1} P_i A$.

\bigskip

We continue in the same setting and show that $g A^{-1} P_i A g^{-1}$ is in fact equal to $P_i$ for $i =1,2$. This proves that $g A^{-1}$ is the identity and so finally every element of $\Gamma$ is already contained in $G$.

We show this by excluding that $P_1$ is conjugated to a proper power of itself, or to a power of $P_2$. First, the cylinder decomposition to which $P_1$ corresponds provides a longest saddle connection in the eigen direction of $P_1$. If $P_1$ would be conjugated to a proper positive power of itself by an element of $\Gamma$ then the affine homeomorphisms of the Chamanara surface corresponding to that element would have to stretch all cylinders and saddle connections in that direction by a factor greater than $1$, contradicting the existence of a longest saddle connection.
Also, $P_1$ cannot be conjugated to $P_1^{-1}$ because they are not even conjugated in the whole of $\PSL_2(\mathbb{R})$.

By the same argument $P_2$ cannot be conjugated in $\Gamma$ to a proper power of itself because there exist two longest saddle connections in the eigen direction of $P_2$.

Finally, $P_1$ and $P_2$ are not conjugated in $\Gamma$ as the boundary of every cylinder in the first cylinder decomposition consists of four saddle connections (for the largest cylinder, one saddle connection is counted with 
multiplicity $2$), whereas the boundary of the largest cylinder in the second cylinder decomposition consists of two saddle connections.

It follows that $g A^{-1} P_i A g^{-1} = P_i$ for $i=1,2$ and so $g A^{-1}$ commutes with both $P_1$ and $P_2$. Since nontrivial Möbius transformations commute if and only if they have the same set of fixed points we conclude that $g A^{-1}$ is the identity.

Hence, we have proven that the Veech group of the Chamanara surface is $G$. We summarize this in a last proposition.

\begin{prop}
 The Veech group $\Gamma$ of the Chamanara surface is generated by the two parabolic elements
 $P_1 =
 \begin{pmatrix}
  1 & 6 \\
  0 & 1
 \end{pmatrix}$ and
 $P_2 =
 \frac{1}{4} \cdot
 \begin{pmatrix}
  -5 & 27 \\
  -3 & 13
 \end{pmatrix}$.
\end{prop}

\bibliographystyle{amsalpha}
\bibliography{/home/anja/Dokumente/Mathematik/literature/BibTex/Literatur}

\providecommand{\bysame}{\leavevmode\hbox to3em{\hrulefill}\thinspace}
\providecommand{\MR}{\relax\ifhmode\unskip\space\fi MR }
\providecommand{\MRhref}[2]{%
  \href{http://www.ams.org/mathscinet-getitem?mr=#1}{#2}
}
\providecommand{\href}[2]{#2}
\begin{thebibliography}{Mas88}

\bibitem[Bea95]{beardon_95}
Alan~F. Beardon, \emph{The geometry of discrete groups}, corr. 2. print. ed.,
  Graduate texts in mathematics; 91, Springer, New York, 1995.

\bibitem[Cha04]{chamanara_04}
Reza Chamanara, \emph{Affine automorphism groups of surfaces of infinite type},
  In the Tradition of Ahlfors and Bers, III (William Abikoff and Andrew Haas,
  eds.), Contemporary mathematics, vol. 355, 2004, pp.~123--145.

\bibitem[Mas88]{maskit_88}
Bernard Maskit, \emph{Kleinian groups}, Die Grundlehren der mathematischen
  Wissenschaften in Einzeldarstellungen; 287, Springer, Berlin, 1988.

\bibitem[Vee89]{veech_89}
William~A. Veech, \emph{Teichmüller curves in moduli space, {Eisenstein}
  series and an application to triangular billiards}, Inventiones Mathematicae
  \textbf{97} (1989), no.~3, 553--583.

\end{thebibliography}

\end{document}